\documentclass[11pt]{article}

\usepackage{latexsym}
\usepackage[dvips]{graphicx}
\usepackage{float}
\usepackage{amsmath}
\usepackage{amssymb}
\usepackage{epic}
\usepackage{color}
\usepackage{setspace}
\usepackage{lscape}
\usepackage{arydshln}
\usepackage{fancybox}   
\usepackage{mathtools}
\usepackage{amsthm}

\begin{document}
\bibliographystyle{plain}
\floatplacement{table}{H}
\newtheorem{definition}{Definition}[section]
\newtheorem{lemma}{Lemma}[section]
\newtheorem{theorem}{Theorem}[section]
\newtheorem{corollary}{Corollary}[section]
\newtheorem{proposition}{Proposition}[section]

\newcommand{\sni}{\sum_{i=1}^{n}}
\newcommand{\snj}{\sum_{j=1}^{n}}
\newcommand{\smj}{\sum_{j=1}^{m}}
\newcommand{\sumjm}{\sum_{j=1}^{m}}
\newcommand{\bdis}{\begin{displaymath}}
\newcommand{\edis}{\end{displaymath}}
\newcommand{\beq}{\begin{equation}}
\newcommand{\eeq}{\end{equation}}
\newcommand{\beqn}{\begin{eqnarray}}
\newcommand{\eeqn}{\end{eqnarray}}
\newcommand{\defeq}{\stackrel{\triangle}{=}}
\newcommand{\simleq}{\stackrel{<}{\sim}}
\newcommand{\sep}{\;\;\;\;\;\; ; \;\;\;\;\;\;}
\newcommand{\real}{\mbox{$ I \hskip -4.0pt R $}}
\newcommand{\complex}{\mbox{$ I \hskip -6.8pt C $}}
\newcommand{\integ}{\mbox{$ Z $}}
\newcommand{\realn}{\real ^{n}}
\newcommand{\sqrn}{\sqrt{n}}
\newcommand{\sqrtwo}{\sqrt{2}}
\newcommand{\prf}{{\bf Proof. }}

\newcommand{\onehlf}{\frac{1}{2}}
\newcommand{\thrhlf}{\frac{3}{2}}
\newcommand{\fivhlf}{\frac{5}{2}}
\newcommand{\onethd}{\frac{1}{3}}
\newcommand{\lb}{\left ( }
\newcommand{\lcb}{\left \{ }
\newcommand{\lsb}{\left [ }
\newcommand{\labs}{\left | }
\newcommand{\rb}{\right ) }
\newcommand{\rcb}{\right \} }
\newcommand{\rsb}{\right ] }
\newcommand{\rabs}{\right | }
\newcommand{\lnm}{\left \| }
\newcommand{\rnm}{\right \| }
\newcommand{\lambdab}{\bar{\lambda}}
%
%
\newcommand{\xj}{x_{j}}
\newcommand{\xjb}{\bar{x}_{j}}
\newcommand{\xro}{x_{\resh}}
\newcommand{\xrob}{\bar{x}_{\resh}}
\newcommand{\xsig}{x_{\sigma}}
\newcommand{\xsigb}{\bar{x}_{\sigma}}
\newcommand{\xnmjb}{\bar{x}_{n-j+1}}
\newcommand{\xnmj}{x_{n-j+1}}
\newcommand{\aroj}{a_{\resh j}}
\newcommand{\arojb}{\bar{a}_{\resh j}}
\newcommand{\aroro}{a_{\resh \resh}}
\newcommand{\amuro}{a_{\mu \resh}}
\newcommand{\amumu}{a_{\mu \mu}}
\newcommand{\aii}{a_{ii}}
\newcommand{\aik}{a_{ik}}
\newcommand{\akj}{a_{kj}}
\newcommand{\atwoii}{a^{(2)}_{ii}}
\newcommand{\atwoij}{a^{(2)}_{ij}}
\newcommand{\ajj}{a_{jj}}
\newcommand{\aiib}{\bar{a}_{ii}}
\newcommand{\ajjb}{\bar{a}_{jj}}
\newcommand{\bii}{a_{jj}}
\newcommand{\biib}{\bar{a}_{jj}}
\newcommand{\aij}{a_{i,n-i+1}}
\newcommand{\akl}{a_{j,n-j+1}}
\newcommand{\aijb}{\bar{a}_{i,n-i+1}}
\newcommand{\aklb}{\bar{a}_{j,n-j+1}}
\newcommand{\bij}{a_{n-j+1,j}}
\newcommand{\arorob}{\bar{a}_{\resh \resh}}
\newcommand{\arosig}{a_{\resh \sigma}}
\newcommand{\arosigb}{\bar{a}_{\resh \sigma}}
\newcommand{\sumjrosig}{\sum_{\stackrel{j=1}{j\neq\resh,\sigma}}^{n}}
\newcommand{\summuro}{\sum_{\stackrel{j=1}{j\neq\mu,\resh}}^{n}}
\newcommand{\sumjnoti}{\sum_{\stackrel{j=1}{j\neq i}}^{n}}
\newcommand{\sumlnoti}{\sum_{\stackrel{\ell=1}{\ell \neq i}}^{n}}
\newcommand{\sumknoti}{\sum_{\stackrel{k=1}{k\neq i}}^{n}}
\newcommand{\sumknotij}{\sum_{\stackrel{k=1}{k\neq i,j}}^{n}}
\newcommand{\sumk}{\sum_{k=1}^{n}}
\newcommand{\snl}{\sum_{\ell=1}^{n}}
\newcommand{\sumji}{\sum_{\stackrel{j=1}{j\neq i, n-i+1}}^{n}}
\newcommand{\sumki}{\sum_{\stackrel{k=1}{k\neq i, n-i+1}}^{n}}
\newcommand{\sumkj}{\sum_{\stackrel{k=1}{k\neq j, n-j+1}}^{n}}
\newcommand{\sumjro}{\sum_{\stackrel{j=1}{j\neq\resh}}^{n}}
\newcommand{\rrosig}{R''_{\resh \sigma}}
\newcommand{\rro}{R'_{\resh}}
\newcommand{\gamror}{\Gamma_{\resh}^{R}(A)}
\newcommand{\gamir}{\Gamma_{i}^{R}(A)}
\newcommand{\gamctrr}{\Gamma_{\frac{n+1}{2}}^{R}(A)}
\newcommand{\gamctrc}{\Gamma_{\frac{n+1}{2}}^{C}(A)}
\newcommand{\gamroc}{\Gamma_{\resh}^{C}(A)}
\newcommand{\gamjc}{\Gamma_{j}^{C}(A)}
\newcommand{\lamror}{\Lambda_{\resh}^{R}(A)}
\newcommand{\lamir}{\Lambda_{i}^{R}(A)}
\newcommand{\lamirepsilon}{\Lambda_{i}^{R}(A_{\epsilon})}
\newcommand{\lamnir}{\Lambda_{n-i+1}^{R}(A)}
\newcommand{\lamjr}{\Lambda_{j}^{R}(A)}
\newcommand{\varphiij}{\Phi_{ij}^{R}(A)}
\newcommand{\delir}{\Delta_{i}^{R}(A)}
\newcommand{\vir}{V_{i}^{R}(A)}
\newcommand{\pamir}{\Pi_{i}^{R}(A)}
\newcommand{\xir}{\Xi_{i}^{R}(A)}
\newcommand{\lamjc}{\Lambda_{j}^{C}(A)}
\newcommand{\vjc}{V_{j}^{C}(A)}
\newcommand{\pamjc}{\Pi_{j}^{C}(A)}
\newcommand{\xjc}{\Xi_{j}^{C}(A)}
\newcommand{\lamroc}{\Lambda_{\resh}^{C}(A)}
\newcommand{\lamsigr}{\Lambda_{\sigma}^{R}(A)}
\newcommand{\lamsigc}{\Lambda_{\sigma}^{C}(A)}
\newcommand{\psii}{\Psi_{i}^{R}(A)}
\newcommand{\psiq}{\Psi_{q}}
\newcommand{\psiiepsilon}{\Psi_{i}^{R}(A_{\epsilon})}
\newcommand{\psiqepsilon}{\Psi_{q}(A_{\epsilon})}
\newcommand{\psiqc}{\Psi_{q}^{c}}
\newcommand{\psiqcepsilon}{\Psi_{q}^{c}(A_{\epsilon})}

\newcommand{\xmu}{x_{\mu}}
\newcommand{\xmub}{\bar{x}_{\mu}}
\newcommand{\xnu}{x_{\nu}}
\newcommand{\xnub}{\bar{x}_{\nu}}
\newcommand{\amuj}{a_{\mu j}}
\newcommand{\amujb}{\bar{a}_{\mu j}}
\newcommand{\amumub}{\bar{a}_{\mu \mu}}
\newcommand{\amunu}{a_{\mu \nu}}
\newcommand{\amunub}{\bar{a}_{\mu \nu}}
\newcommand{\sumjmunu}{\sum_{\stackrel{j=1}{j\neq\mu,\nu}}}
\newcommand{\rmunu}{R''_{\mu \nu}}
\newcommand{\rmu}{R'_{\mu}}

\newcommand{\Azero}{A_{0}}
\newcommand{\Aone}{A_{1}}
\newcommand{\Atwo}{A_{2}}
\newcommand{\Ath}{A_{3}}
\newcommand{\Afr}{A_{4}}
\newcommand{\Afv}{A_{5}}
\newcommand{\Asx}{A_{6}}
\newcommand{\Anmo}{A_{n-1}}
\newcommand{\Anmt}{A_{n-2}}
\newcommand{\An}{A_{n}}
\newcommand{\Aj}{A_{j}}
\newcommand{\Anmell}{A_{n-\ell}}
\newcommand{\Anmellmo}{A_{n-\ell-1}}
\newcommand{\Anmkmo}{A_{n-k-1}}
\newcommand{\Anmkmell}{A_{n-k-\ell}}
\newcommand{\Anmtk}{A_{n-2k}}
\newcommand{\Anmkmellmo}{A_{n-k-\ell-1}}

\newcommand{\anmell}{a_{n-\ell}}
\newcommand{\anmellmo}{a_{n-\ell-1}}
\newcommand{\anmkmo}{a_{n-k-1}}
\newcommand{\anmkmell}{a_{n-k-\ell}}
\newcommand{\anmkmellmo}{a_{n-k-\ell-1}}

\newcommand{\azero}{a_{0}}
\newcommand{\aone}{a_{1}}
\newcommand{\atwo}{a_{2}}
\newcommand{\ath}{a_{3}}
\newcommand{\afr}{a_{4}}
\newcommand{\afv}{a_{5}}
\newcommand{\asx}{a_{6}}
\newcommand{\anmo}{a_{n-1}}
\newcommand{\anmt}{a_{n-2}}
\newcommand{\an}{a_{n}}
\newcommand{\aj}{a_{j}}
\newcommand{\ajmo}{a_{j-1}}
\newcommand{\ajmt}{a_{j-2}}

\newcommand{\Bzero}{B_{0}}
\newcommand{\Bone}{B_{1}}
\newcommand{\Btwo}{B_{2}}
\newcommand{\Bth}{B_{3}}
\newcommand{\Bfr}{B_{4}}
\newcommand{\Bfv}{B_{5}}
\newcommand{\Bsx}{B_{6}}
\newcommand{\Bnmo}{B_{n-1}}
\newcommand{\Bndtwomo}{B_{n/2-1}}
\newcommand{\Bnmt}{B_{n-2}}
\newcommand{\Bn}{B_{n}}
\newcommand{\Bj}{B_{j}}

\newcommand{\Anmk}{A_{n-k}}
\newcommand{\Bnmk}{A_{n-k}}
\newcommand{\anmk}{a_{n-k}}
\newcommand{\zj}{z^{j}}          
\newcommand{\zk}{z^{k}}          
\newcommand{\zkpm}{z^{k+m}}          
\newcommand{\zm}{z^{m}}          

\newcommand{\cii}{c_{ii}}
\newcommand{\cik}{c_{ik}}
\newcommand{\ckj}{c_{kj}}
\newcommand{\ctwoii}{c^{(2)}_{ii}}
\newcommand{\ctwoij}{c^{(2)}_{ij}}
\newcommand{\cjj}{c_{jj}}

\newcommand{\bik}{b_{ik}}
\newcommand{\bkj}{b_{kj}}
\newcommand{\btwoii}{b^{(2)}_{ii}}
\newcommand{\btwoij}{b^{(2)}_{ij}}
\newcommand{\bjj}{b_{jj}}

\newcommand{\abii}{(AB)_{ii}}
\newcommand{\abil}{(AB)_{i\ell}}

\newcommand{\bkl}{b_{k\ell}}
\newcommand{\btwoil}{b^{(2)}_{i\ell}}
\newcommand{\bll}{b_{\ell \ell}}

\newcommand{\matrixspace}{\;\;}
\newcommand{\ellone}{\ell_{1}}
\newcommand{\elltwo}{\ell_{2}}

\newcommand{\varphik}{\varphi_{k}}
\newcommand{\chik}{\chi_{k}}
\newcommand{\Phik}{\Phi_{k}}
\newcommand{\psik}{\psi_{k}}
\newcommand{\dr}{\beta^{n-k}}
\newcommand{\dn}{\beta^{-k}}
\newcommand{\betatwotilde}{\tilde{\rtwonk}}
\newcommand{\betathtilde}{\tilde{\ronenk}}
\newcommand{\betamink}{\beta^{-k}}

\newcommand{\rvarphik}{r(\varphi_{k})}    
\newcommand{\rrvarphik}{(r(\varphi_{k}))}    
\newcommand{\rchik}{r(\chi_{k})}    
\newcommand{\rrchik}{(r(\chi_{k}))}    
\newcommand{\svarphik}{s(\varphi_{k})}    
\newcommand{\ssvarphik}{(s(\varphi_{k}))}    
\newcommand{\schik}{s(\chi_{k})}    
\newcommand{\sschik}{(s(\chi_{k}))}    

\newcommand{\rpsik}{r(\psi_{k})}    
\newcommand{\rrpsik}{(r(\psi_{k}))} 
\newcommand{\spsik}{s(\psi_{k})}    
\newcommand{\sspsik}{(s(\psi_{k}))} 
\newcommand{\rpsinmk}{r(\psi_{n-k})}    
\newcommand{\rrpsinmk}{(r(\psi_{n-k}))} 
\newcommand{\spsinmk}{s(\psi_{n-k})}    
\newcommand{\sspsinmk}{(s(\psi_{n-k}))} 

\newcommand{\rvarphinmk}{r(\varphi_{n-k})}    
\newcommand{\rrvarphinmk}{(r(\varphi_{n-k}))}    
\newcommand{\rchinmk}{r(\chi_{n-k})}    
\newcommand{\rrchinmk}{(r(\chi_{n-k}))}    
\newcommand{\svarphinmk}{s(\varphi_{n-k})}    
\newcommand{\ssvarphinmk}{(s(\varphi_{n-k}))}    
\newcommand{\schinmk}{s(\chi_{n-k})}    
\newcommand{\sschinmk}{(s(\chi_{n-k}))}    
\newcommand{\stilde}{\tilde{s}}

\newcommand{\resh}{\rho}
\newcommand{\snk}{s}
\newcommand{\ronenk}{r_{1}}
\newcommand{\rtwonk}{r_{2}}
\newcommand{\obar}{\widebar{O}}
\newcommand{\Lbar}{\widebar{L}}
\newcommand{\dltaone}{\delta_{1}}
\newcommand{\dltatwo}{\delta_{2}}
\newcommand{\dltatld}{\tilde{\delta}}
\newcommand{\tauone}{\tau_{1}}
\newcommand{\tautwo}{\tau_{2}}
\newcommand{\xb}{\bar{x}}
\newcommand{\qone}{q_{1}}
\newcommand{\qtwo}{q_{2}}
\newcommand{\ub}{\bar{u}}
\newcommand{\cmm}{\complex^{m \times m}}
\newcommand{\opdsk}{\mathcal{O}}
\newcommand{\cldsk}{\widebar{\mathcal{O}}}

\newcommand{\sqrtaone}{\sqrt{\a1}}
\newcommand{\sqrtatwo}{\sqrt{\atwo}}

\newcommand{\rhoone}{\rho_{1}}
\newcommand{\rhotwo}{\rho_{2}}
\newcommand{\xone}{x_{1}}
\newcommand{\xtwo}{x_{2}}
\newcommand{\xthr}{x_{3}}
\newcommand{\xfor}{x_{4}}
\newcommand{\xfiv}{x_{5}}

\newcommand{\lnorm}{\left \|}
\newcommand{\rnorm}{\right \|}
\newcommand{\lnrm}{\biggl | \biggl |}
\newcommand{\rnrm}{\biggr |\biggr |}

\newcommand{\recipp}{p^{\protect \#}}
\newcommand{\recipP}{P^{\protect \#}}

\newcommand{\absatwo}{|\atwo|}
\newcommand{\absqrtatwo}{\sqrt{|\atwo|}}
\newcommand{\absrhoone}{|\rhoone|}
\newcommand{\absrhotwo}{|\rhotwo|}
\newcommand{\sqrgamma}{\sqrt{\gamma}}
\newcommand{\aminb}{a-b}                 
\newcommand{\sqrtc}{\sqrt{c}}            
\newcommand{\sqrtabsc}{\sqrt{|c|}}            
\newcommand{\absalf}{|\alpha|}            
\newcommand{\rhoi}{\rho_{i}}                  
\newcommand{\sigmai}{\sigma_{i}}                  
\newcommand{\rhoj}{\rho_{j}}                  
\newcommand{\sigmaj}{\sigma_{j}}                  
\newcommand{\xij}{x_{ij}}                  
\newcommand{\yij}{y_{ij}}                  
\newcommand{\varphii}{\varphi_{i}}
\newcommand{\varphij}{\varphi_{j}}
\newcommand{\varphijmo}{\varphi_{j-1}}
\newcommand{\varphijmt}{\varphi_{j-2}}
\newcommand{\varphin}{\varphi_{n}}
\newcommand{\varphinmo}{\varphi_{n-1}}
\newcommand{\varphinmt}{\varphi_{n-2}}
\newcommand{\varphione}{\varphi_{1}}
\newcommand{\varphizero}{\varphi_{0}}
\newcommand{\varphitwo}{\varphi_{2}}
\newcommand{\psione}{\psi_{1}}
\newcommand{\psitwo}{\psi_{2}}
\newcommand{\xstar}{x^{\ast}}
\newcommand{\uuone}{\mathcal{U}_{1}}
\newcommand{\uutwo}{\mathcal{U}_{2}}
\newcommand{\xkova}{\hat{x}}
\newcommand{\xkovalam}{\xkova_{max}(\ell-1)}

\begin{center}
\large
{\bf IMPROVED CAUCHY RADIUS FOR SCALAR AND MATRIX POLYNOMIALS}
\vskip 0.5cm
\normalsize
A. Melman \\
Department of Applied Mathematics \\
School of Engineering, Santa Clara University  \\
Santa Clara, CA 95053  \\
e-mail : amelman@scu.edu \\
\vskip 0.5cm
\end{center}

\begin{abstract}
We improve the Cauchy radius of both scalar and matrix polynomials, which is an upper bound on the moduli
of the zeros and eigenvalues, respectively, by using appropriate polynomial multipliers.
\vskip 0.15cm
{\bf AMS(MOS) subject classification :} 30C15, 47A56, 65F15
\end{abstract}

%
%
%
%

\section{Introduction}           
\label{introduction} 

A simple but classical result from~1829 due to Cauchy (\cite{Cauchy}, \cite[Th.(27,1), p.122 and Exercise 1, p.126]{Marden})
states that the zeros of a polynomial $p(z)= a_{n}z^{n} + a_{n-1} z^{n-1} + \dots + a_{1} z + a_{0}$, with complex coefficients and $\an \neq 0$,
lie in $|z| \leq \rho[p]$, 
where $\rho[p]$ is the \emph{Cauchy radius} of $p$, namely, the unique positive solution of 
\bdis
|a_{n}| z^{n} - |a_{n-1}| z^{n-1} - \dots  - |\aone| z - |\azero| = 0  \; . 
\edis
A smaller Cauchy radius was obtained much more recently by Rahman and Schmeisser (\cite[Theorem~8.3.1]{RS}), who showed that
$\rho[(\an z^{k}- \anmk)p(z)] \leq \rho[p]$, where $k$ is the smallest positive integer such that $\anmk \neq 0$, i.e., a better bound can be found by 
using a polynomial multiplier. 

A generalization to matrix polynomials of Cauchy's classical bound for scalar polynomials
was derived in~\cite{BiniNoferiniSharify}, \cite{HighamTisseur}, and \cite{Melman_MatPol}. It states that 
all the eigenvalues of the regular matrix polynomial $P(z) = \An z^{n} + A_{n-1}z^{n-1} + \dots + A_{1}z + A_{0}$, with complex coefficient matrices
and $\An$ nonsingular, lie in $|z| \leq \rho[P]$, where, as in the scalar case, $\rho[P]$ is called the \emph{Cauchy radius} of $P$, 
which is the unique positive solution of
\bdis
\|\An^{-1} \|^{-1} z^{n} - \|A_{n-1}\|z^{n-1} - \dots - \|A_{1}\|z - \|A_{0}\| = 0
\edis
for any matrix norm. The eigenvalues of~$P$ are the complex numbers $z$ for which a nonzero complex vector $v$ exists such that $P(z)v=0$.
If $\An$ is nonsingular, they are the solutions of $\emph{det}{P(z)}=0$. A matrix polynomial $P$ is \emph{regular} if $\emph{det}{P}$ is not identically zero.
When $P$ is linear and monic, i.e., $P(z)=Iz - \Azero$, one obtains the standard eigenvalue problem.

In~\cite{Melman_RS}, the improved Cauchy radius from~\cite{RS} was also generalized to matrix polynomials. It was shown there that, under mild 
conditions on $\An$, both $\rho \lsb \lb \An z^{k} - \Anmk \rb P(z) \rsb \leq \rho \lsb P \rsb$ and  
$\rho \lsb P(z) \lb \An z^{k} - \Anmk \rb \rsb \leq \rho \lsb P \rsb$, with $k$ the smallest positive integer such that $\Anmk$ is not the null matrix.

There do not seem to exist other multipliers with these properties in the literature, and our purpose here is to derive different 
multipliers that also improve the Cauchy radius for both scalar and matrix polynomials and that, in general, perform better than the improvements 
from~\cite{Melman_RS} and~\cite{RS}.

In Section~\ref{matrixpolynomials} we present such polynomial multipliers first for matrix
polynomials, while we consider scalar polynomials as a special case in Section~\ref{scalarpolynomials}.

%
%
%
%

\section{Improved Cauchy radius for matrix polynomials}
\label{matrixpolynomials}

The following theorem presents three matrix polynomials, obtained by multiplying a given matrix polynomial $P$ by another matrix polynomial, that have a 
smaller Cauchy radius than that of $P$. Clearly, for any matrix polynomial $T$, a region in the complex plane containing all the eigenvalues of $TP$ or $PT$ 
also contains those of $P$.

%
%
\begin{theorem}
\label{Theorem_Matrix_Bakchik}
Let $P(z) = \sum_{j=0}^{n} A_{j} z^{j}$ be a regular matrix polynomial of degree $n$ that is at least a trinomial, with square complex matrices 
$A_{j}$ ($0 \leq j \leq n$) and $\An$ nonsingular, and let~$k$~and~$\ell$ be the smallest positive integers such that $\Anmk$ and
$\Anmkmell$ are not the null matrix. Define
\begin{eqnarray*}
& & Q^{(L)}_{1}(z) = \lb \An z^{k+\ell} - \Anmk z^{\ell} - \Anmkmell \rb P(z) \; ,    \\
& & Q^{(L)}_{2}(z) = \lb \An z^{2k} - \Anmk z^{k} + \Anmk^{2}\An^{-1} \rb P(z) \; ,   \\
\hskip -1cm \text{and, when $\ell=k$,} & &   \\ 
& & Q^{(L)}_{3}(z) = \lb \An z^{2k} - \Anmk z^{k} - A_{n-2k} + \Anmk^{2}\An^{-1} \rb P(z) \; .     
\end{eqnarray*}
Furthermore, define 
\begin{eqnarray*}
& & Q^{(R)}_{1}(z) = P(z) \lb \An z^{k+\ell} - \Anmk z^{\ell} - \Anmkmell \rb \; ,     \\
& & Q^{(R)}_{2}(z) = P(z) \lb \An z^{2k} - \Anmk z^{k} + \Anmk^{2} \An^{-1} \rb  \; ,    \\
\hskip -1cm \text{and, when $\ell=k$,} & &   \\ 
& & Q^{(R)}_{3}(z) = P(z) \lb \An z^{2k} - \Anmk z^{k} - A_{n-2k} + \Anmk^{2} \An^{-1} \rb \; . 
\end{eqnarray*}
\noindent
For any matrix norm $\|.\|$, if $\|\An^{-2}\|^{-1} = \|\An\|\|\An^{-1}\|^{-1}$ and if 
$\An\Anmk=\Anmk\An$, $\An\Anmkmell = \Anmkmell \An$,
then, for any admissible values of $k$ and $\ell$, it follows that $\rho[Q^{(L)}_{j}] \leq \rho[P]\,$ for $j=1,2$.
If $\ell=k$, then $\rho[Q^{(L)}_{3}] \leq \rho[P]$. Analogous results hold for $Q^{(R)}_{j}$ ($j=1,2,3$).
\end{theorem}
\prf
We prove the theorem for $Q_{j}^{(L)}$ ($j=1,2,3$); the proof for $Q_{j}^{(R)}$ ($j=1,2,3$) is analogous.
Let $k$ and $\ell$ be as in the statement of the theorem and, if it exists, let $s$ be the first positive integer such that $A_{n-k-\ell-s}$ is 
not the null matrix (when $P$ is a trinomial, then no such $s$ exists). We will make use of the following expression, where $m$ is any positive 
integer and $M$ is any square complex matrix:
\begin{eqnarray}
\label{multexpression}
& & \hskip -1cm \lb \An z^{k+m} - \Anmk z^{m} - M \rb P(z)  \nonumber \\
& & = \An^{2} z^{n+k+m} + \An\Anmkmell z^{n-\ell+m} -\Anmk^{2} z^{n-k+m} - M \An z^{n} \nonumber \\
& & \hskip 4cm + \An z^{k+m} \sum_{j=0}^{n-k-\ell-s} \Aj z^{j}
- \Anmk z^{m} \sum_{j=0}^{n-k-\ell} A_{j} z^{j}
- M \sum_{j=0}^{n-k} \Aj z^{j}  \; . \nonumber \\ 
\end{eqnarray}
If $P$ is a trinomial, then the summation in~(\ref{multexpression}) with upper index limit $n-k-\ell-s$ is set equal to zero.

We begin with $Q^{(L)}_{1}$, which is obtained by setting $m=\ell$ and $M=\Anmkmell$ in~(\ref{multexpression}):
\beq
\label{powereq1}
\lb \An z^{k+\ell} - \Anmk z^{\ell} - \Anmkmell \rb P(z)  = Q^{(L)}_{1}(z) = \An^{2} z^{n+k+\ell} + S(z) \; ,
\eeq
where 
\beq
\label{powereq2}
S(z) = -\Anmk^{2} z^{n-k+\ell} + \An z^{k+\ell} \sum_{j=0}^{n-k-\ell-s} \Aj z^{j}
- \Anmk z^{\ell} \sum_{j=0}^{n-k-\ell} A_{j} z^{j} - \Anmkmell \sum_{j=0}^{n-k} \Aj z^{j} = \sum_{j=0}^{\nu} B_{j}z^{j} \; ,
\eeq
$\nu \leq n-\min\{s,k-\ell\}$ (when $P$ is a trinomial, $\nu \leq 2(n-k)$, since in this case $k+\ell=n$ so that 
$k-\ell=2k-n$), and each matrix $B_{j}$ is a sum of terms of the form $A_{j}$ or $A_{i}A_{j}$.
If we define
\bdis
\Phi(z) = \sum_{j=0}^{\nu} \|B_{j}\|z^{j} \; ,
\edis
then the Cauchy radius of $Q^{(L)}_{1}$ is the unique positive solution of $\|\An^{-2}\|^{-1} z^{n+k+\ell} - \Phi(z) = 0$. 
We now set $x=\rho[P]$, i.e., $x$ satisfies 
\beq
\label{rhoP}
\|\An^{-1}\|^{-1} x^{n} - \|\Anmk\|x^{n-k} - \|\Anmkmell\| x^{n-k-\ell} - \sum_{j=0}^{n-k-\ell-s} \|\Aj\|x^{j} = 0 \; .
\eeq
Using~(\ref{rhoP}) and the basic properties $\|A+B\| \leq \|A\|+\|B\|$ and $\|AB\| \leq \|A\|\|B\|$ of matrix norms, we have that
\begin{eqnarray}
\Phi(x) & \leq & \|\Anmk\|^{2} x^{n-k+\ell} + \|\An\| x^{k+\ell} \sum_{j=0}^{n-k-\ell-s} \|\Aj\| x^{j} \label{firstsummation}   \\ 
        &      &  \hskip 4cm + \|\Anmk\| x^{\ell} \sum_{j=0}^{n-k-\ell} \|A_{j}\| x^{j} + \|\Anmkmell\| \sum_{j=0}^{n-k} \|\Aj\| x^{j} \nonumber \\
        &   =  & \|\Anmk\|^{2} x^{n-k+\ell} + \|\An\| x^{k+\ell} \lb \|\An^{-1}\|^{-1}x^{n} -\|\Anmk\|x^{n-k} - \|\Anmkmell\| x^{n-k-\ell} \rb  \nonumber \\[8pt]
        &      & \hskip 2cm + \|\Anmk\| x^{\ell} \lb \|\An^{-1}\|^{-1} x^{n} -\|\Anmk\|x^{n-k} \rb + \|\Anmkmell\| \lb \|\An^{-1}\|^{-1} x^{n} \rb  \nonumber \\[8pt]
        & \leq &  \|\An\|\|\An^{-1}\|^{-1} x^{n+k+\ell}  = \|\An^{-2}\|^{-1} x^{n+k+\ell}  \; . \nonumber 
\end{eqnarray}
We have used both the fact that $\|\An^{-1}\|^{-1} \leq \|\An\|$ and our assumption that $\|\An^{-2}\|^{-1} =  \|\An\|\|\An^{-1}\|^{-1}$.
This means that $\|\An^{-2}\|^{-1} x^{n+k+\ell} - \Phi(x) \geq 0$ and, therefore, that $x$ must lie to the right of $\rho[Q^{(L)}_{1}]$,
i.e., $\rho[Q^{(L)}_{1}] \leq \rho[P]$. 

When $P$ is a trinomial,
the second term in the right-hand side of~(\ref{firstsummation}) is absent, and the result follows analogously. 


For $Q^{(L)}_{2}$ the proof is similar and we will omit unnecessary details. Here we set $m=k$ and $M=-\Anmk^{2}\An^{-1}$ 
in~(\ref{multexpression}) to obtain 
\beq
\label{powereq3}
\lb \An z^{2k} - \Anmk z^{k} + \Anmk^{2} \An^{-1} \rb P(z) = Q^{(L)}_{2}(z) = \An^{2} z^{n+2k} + S(z) \; ,
\eeq
where
\beq
\label{powereq4}
S(z) = \An\Anmkmell z^{n-\ell+k} + \An z^{2k} \sum_{j=0}^{n-k-\ell-s} \Aj z^{j}
- \Anmk z^{k} \sum_{j=0}^{n-k-\ell} A_{j} z^{j} + \Anmk^{2}\An^{-1} \sum_{j=0}^{n-k} \Aj z^{j} = \sum_{j=0}^{\nu} B_{j}z^{j} \; ,
\eeq
and $\nu = n-\min\{k, \ell-k \}$ (when $P$ is a trinomial, $\nu = n-\min\{k, n-2k\}$).
The Cauchy radius of $Q^{(L)}_{2}$ is the unique positive solution of $\|\An^{-2}\|^{-1} z^{n+2k} - \Phi(z) = 0$,
where $\Phi(z) = \sum_{j=0}^{\nu} \|B_{j}\|z^{j}$.
With $x=\rho[P]$ and~(\ref{rhoP}), we have 
\begin{eqnarray}
\Phi(x) & \leq & \|\An\|\|\Anmkmell\| x^{n-\ell+k} + \|\An\| x^{2k} \sum_{j=0}^{n-k-\ell-s} \|\Aj\| x^{j}  \nonumber \\
        &      & \hskip 4cm + \|\Anmk\| x^{k} \sum_{j=0}^{n-k-\ell} \|A_{j}\| x^{j} + \|\Anmk\|^{2}\|\An^{-1}\| \sum_{j=0}^{n-k} \|\Aj\| x^{j} \nonumber \\[8pt]
        &   =  & \|\An\|\|\Anmkmell\| x^{n-\ell+k} + \|\An\| x^{2k} \lb \|\An^{-1}\|^{-1}x^{n} -\|\Anmk\|x^{n-k} - \|\Anmkmell\| x^{n-k-\ell} \rb  \nonumber \\[8pt]
        &      & \hskip 2cm + \|\Anmk\| x^{k} \lb \|\An^{-1}\|^{-1} x^{n} -\|\Anmk\|x^{n-k} \rb 
                 + \|\Anmk\|^{2}\|\An^{-1}\| \lb \|\An^{-1}\|^{-1} x^{n} \rb \nonumber \\[8pt]
        & \leq &  \|\An\|\|\An^{-1}\|^{-1} x^{n+2k}  = \|\An^{-2}\|^{-1} x^{n+2k}  \; . \nonumber 
\end{eqnarray}
Therefore, $\|\An^{-2}\|^{-1} x^{n+2k} - \Phi(x) \geq 0$, implying that the Cauchy radius of $Q^{(L)}_{2}$ is smaller than that of $P$. 
When $P$ is a trinomial, the same result follows analogously as before.

For $Q^{(L)}_{3}$, with $\ell=k$, we set $m=k=\ell$ and $M=A_{n-2k} - \Anmk^{2} \An^{-1}$ in~(\ref{multexpression}), which gives
\beq
\label{powereq5}
\lb \An z^{2k} - \Anmk z^{k} - \Anmtk + \Anmk^{2} \An^{-1} \rb P(z)  = Q^{(L)}_{3}(z) = \An^{2} z^{n+2k} + S(z) \; ,
\eeq
where
\beq
\label{powereq6}
S(z) = \An z^{2k} \sum_{j=0}^{n-2k-s} \Aj z^{j}
- \Anmk z^{k} \sum_{j=0}^{n-2k} A_{j} z^{j} - \lb \Anmtk - \Anmk^{2}\An^{-1} \rb \sum_{j=0}^{n-k} \Aj z^{j} = \sum_{j=0}^{\nu} B_{j}z^{j} \; ,
\eeq
and $\nu \leq n-\min\{k, s \} \leq n-1$ ($\nu \leq n-k \leq n-1$ when $P$ is a trinomial).
With $\Phi(z) = \sum_{j=0}^{\nu} \|B_{j}\|z^{j}$, the Cauchy radius of $Q^{(L)}_{3}$ is the unique positive solution of~$\|\An^{-2}\|^{-1}~z^{n+2k}~-~\Phi(z)~=~0$. 
Setting $x=\rho[P]$, which satisfies equation~(\ref{rhoP}), we have 
\begin{eqnarray}
\Phi(x) & \leq & \|\An\| x^{2k} \sum_{j=0}^{n-2k-s} \|\Aj\| x^{j}
                  + \|\Anmk\| x^{k} \sum_{j=0}^{n-2k} \|A_{j}\| x^{j} + \|\Anmtk - \Anmk^{2}\An^{-1} \| \sum_{j=0}^{n-k} \|\Aj\| x^{j} \nonumber \\
        &   =  & \|\An\| x^{2k} \lb \|\An^{-1}\|^{-1}x^{n} -\|\Anmk\|x^{n-k} - \|\Anmtk\| x^{n-2k} \rb  \nonumber \\[8pt]
        &      &  + \hskip 0.5cm \|\Anmk\| x^{k} \lb \|\An^{-1}\|^{-1} x^{n} -\|\Anmk\|x^{n-k} \rb + \|\Anmtk - \Anmk^{2}\An^{-1} \| \lb \|\An^{-1}\|^{-1} x^{n} \rb  \nonumber \\[8pt]
        & \leq &  \|\An\|\|\An^{-1}\|^{-1} x^{n+2k} + \lb \|\Anmtk - \Anmk^{2}\An^{-1} \| \|\An^{-1}\|^{-1} - \|\An\|\|\Anmtk\| - \|\Anmk\|^{2} \rb x^{n}  \nonumber \\[8pt]
        & \leq &  \|\An\|\|\An^{-1}\|^{-1} x^{n+2k} + \lb \|\Anmtk\|\|\An^{-1}\|^{-1} + \|\Anmk\|^{2}\|\An^{-1}\| \|\An^{-1}\|^{-1}  \right . \nonumber \\[8pt]
        &      &  \hskip 8cm \left . - \|\An\|\|\Anmtk\| - \|\Anmk\|^{2} \rb x^{n}  \nonumber \\[8pt]
        & \leq &  \|\An\|\|\An^{-1}\|^{-1} x^{n+2k} + \lb \|\Anmtk\|\|\An^{-1}\|^{-1} - \|\An\|\|\Anmtk\| \rb x^{n}  \nonumber \\[8pt]
        & \leq &  \|\An\|\|\An^{-1}\|^{-1} x^{n+2k}  = \|\An^{-2}\|^{-1} x^{n+2k}  \; . \nonumber 
\end{eqnarray}
We have obtained that $\|\An^{-2}\|^{-1} x^{n+2k} - \Phi(x) \geq 0$, which means that the Cauchy radius of $Q^{(L)}_{3}$ is smaller than that of $P$. 
This completes the proof. \qed
\vskip 0.25cm

\noindent {\bf Remarks.}
\begin{itemize}
\item 
The matrix $\Anmk^{2}\An^{-1}$ in the definitions of $Q^{(L)}_{j}$ and $Q^{(R)}_{j}$ for $j=2,3$ could be replaced by
$\An^{-1}\Anmk^{2}$ since positive and negative powers (if they exist) of commuting matrices also commute.
\item 
If $A_{n-2k} = \Anmk^{2}\An^{-1}$, then $\rho[Q^{(L)}_{3}]$ and $\rho[Q^{(R)}_{3}]$ are both equal to the improved Cauchy radius of Theorem~8.3.1
in~\cite{RS}.
\item
The conditions $\An\Anmk=\Anmk\An$, $\An\Anmkmell = \Anmkmell \An$, and $\|\An^{-2}\|^{-1} = \|\An\|\|\An^{-1}\|^{-1}$ 
may appear restrictive, but they are always be satisfied if $\An=I$. This can be achieved by multiplying $P$ by $\An^{-1}$, 
which needs to be computed anyway to obtain the Cauchy radius.
\item
In general, the multipliers of $P$ are different from the ones obtained by repeatedly using Theorem~2.2 in~\cite{Melman_RS}, as their degrees can 
easily be seen to be different. 
\end{itemize}

The more zero coefficients a polynomial has, all else being the same, the smaller its Cauchy radius will be. 
Although the matrix polynomials $Q^{(L)}_{j}$ for $j=1,2,3$ may have additional zero coefficients (null matrices),
the ones that we have some control over are the leading zeros immediately following the highest coefficient.
The following lemma allows us to compare their number, thereby indicating which multiplier might be preferable for given values of $k$ and $\ell$.

%
%
\begin{lemma}    
\label{Lemma_powers}                                
Let $P$, $Q^{(L)}_{j}$ and $Q^{(R)}_{j}$ ($j=1,2,3$) be as in Theorem~\ref{Theorem_Matrix_Bakchik}. Then the following holds.
\begin{itemize}         
\item When $\ell < k$, the leading powers of $z$ in $Q^{(L)}_{1}$ are $n+k+\ell$ and $\nu \leq n-1$, whereas for $Q^{(L)}_{2}$ they are
      $n+2k$ and $n+k-\ell \geq n+1$. 
\item When $\ell > k$, the leading powers of $z$ in $Q^{(L)}_{1}$ are $n+k+\ell$ and $n+\ell-k \geq n+1$, whereas for $Q^{(L)}_{2}$ they are
      $n+2k$ and $\nu \leq n-1$.
\item When $\ell = k$, the leading powers of $z$ in $Q^{(L)}_{1}$ and $Q^{(L)}_{2}$ are $n+2k$ and $n$, while for $Q^{(L)}_{3}$, they are $n+2k$ and $\nu \leq n-1$.
\item All of the above results also hold true for $Q^{(R)}_{j}$ ($j=1,2,3$).
\end{itemize}   
\end{lemma}   
\prf
From~(\ref{powereq1}) and~(\ref{powereq2}) we have that, when $\ell < k$, then the leading powers of $Q^{(L)}_{1}$ are $n+k+\ell$ and 
$\nu \leq n-\min\{k-\ell,s\} \leq n-1$ ($\nu \leq n-(k-\ell) \leq n-1$ when $P$ is a trinomial, in which case $k+\ell=n$ so that 
$n-(k-\ell)=2(n-k)$). Here, $s$ is as in the proof of Theorem~\ref{Theorem_Matrix_Bakchik}.
When $\ell > k$, then those powers become $n+k+\ell$ and $n+\ell-k \geq n+1$, and when $\ell=k$, they are
are $n+2k$ and $n$.

Similarly, we observe from~(\ref{powereq3}) and~(\ref{powereq4}) that, when $\ell < k$, the leading powers of $Q^{(L)}_{2}$ are $n+2k$ and $n+k-\ell \geq n+1$,
whereas for $\ell > k$, they are $n+2k$ and $\nu = n-\min\{k, \ell-k \} \leq n-1$. When $\ell = k$, those powers become $n+2k$ and $n$, as for $Q^{(L)}_{1}$.

When $\ell=k$, equations~(\ref{powereq5}) and~(\ref{powereq6}) show that the highest powers of $Q^{(L)}_{3}$ are 
$n+2k$ and $\nu \leq n-\min\{k, s \} \leq n-1$ ($\nu \leq n-k \leq n-1$ when $P$ is a trinomial), with $s$ as in the proof of Theorem~\ref{Theorem_Matrix_Bakchik}.
The proof for $Q^{(R)}_{j}$ ($j=1,2,3$) is analogous. \qed

The number of leading zero coefficients is now easily determined with Lemma~\ref{Lemma_powers} from the leading powers of $Q^{(L)}_{j}$ for $j=1,2,3$.
They can be found on the left in~Table~\ref{Table_zeros} for the worst case (i.e., smallest number of zeros), namely, when $\nu=n-1$, where $\nu$ is 
as in Lemma~\ref{Lemma_powers}, while the degrees of $Q^{(L)}_{j}$ for $j=1,2,3$ can be found on the right.
\begin{table}[H]
\begin{center}
\small
\begin{tabular}{cc}
\hskip -1cm
\begin{tabular}{c|ccc}          
                  & $Q^{(L)}_{1}$   & $Q^{(L)}_{2}$   &  $Q^{(L)}_{3}$  \\        \hline 
                  &                 &                 &                 \\
$\ell < k$        &  $k+\ell$       &  $k+\ell-1$     &     -           \\        
                  &                 &                 &                 \\
$\ell > k$        &  $2k-1$         &  $2k$           &     -           \\        
                  &                 &                 &                 \\
$\ell = k$        &  $2k-1$         &  $2k-1$         &     $2k$        \\        
\end{tabular}
&
\hskip 1cm
\begin{tabular}{c|ccc}          
                  & $Q^{(L)}_{1}$   & $Q^{(L)}_{2}$   &  $Q^{(L)}_{3}$  \\        \hline 
                  &                 &                 &                 \\
$\ell < k$        &  $n+k+\ell$     &  $n+2k$         &     -           \\        
                  &                 &                 &                 \\
$\ell > k$        &  $n+k+\ell$     &  $n+2k$         &     -           \\        
                  &                 &                 &                 \\
$\ell = k$        &  $n+2k$         &  $n+2k$         &     $n+2k$      \\        
\end{tabular}
\\
\end{tabular}
\caption{Number of zero coefficients (left) and degrees of $Q^{(L)}_{j}$ for $j=1,2,3$ (right).}
\label{Table_zeros}     
\end{center}
\end{table}
\normalsize
Table~\ref{Table_zeros} shows that, when $\ell < k$, $Q^{(L)}_{1}$ has a higher number of leading zero coefficients than $Q^{(L)}_{2}$, while its degree
is lower. When $\ell > k$, the same conclusion holds for $Q^{(L)}_{2}$, and when $\ell=k$, then $Q^{(L)}_{3}$ has more such zero coefficients than both 
$Q^{(L)}_{1}$ and $Q^{(L)}_{2}$, while they all have the same degree.
Analogous results
are obtained for $Q^{(R)}_{j}$ for $j=1,2,3$. We thus arrive at the following choice to improve the Cauchy radius of $P$: 
\begin{eqnarray}
& & Q^{(L)}(z) = \lcb 
\begin{array}{ll}
\lb \An z^{k+\ell} - \Anmk z^{\ell} - \Anmkmell \rb P(z)   & \text{if $\ell < k$} \; ,  \\
                                                           &                       \\
\lb \An z^{2k} - \Anmk z^{k} + \Anmk^{2}\An^{-1} \rb P(z)  & \text{if $\ell > k$}  \; , \\
                                                           &                       \\
\lb \An z^{2k} - \Anmk z^{k} - A_{n-2k} + \Anmk^{2}\An^{-1} \rb P(z)  & \text{if $\ell = k$} \; , \\
\end{array}
\right . \nonumber \\ \label{Qdef}
\end{eqnarray}
and we choose $Q^{(R)}$ analogously.

\noindent {\bf Remarks.}
\begin{itemize}

\item
Theorem~\ref{Theorem_Matrix_Bakchik} can be applied recursively to improve the Cauchy radius further. 
One could also alternate between (L) and (R) versions, although, in general, there does not seem to be a large difference between the two.
\item
The improved Cauchy radii require additional matrix multiplications, while a real scalar polynomial equation 
of a degree higher than that of~$P$ needs to be solved. The latter can be dealt with very efficiently so that, as the matrix size increases, the cost
tends to be dominated by the matrix multiplications. It therefore depends on the application if this additional computational cost is justified. 
\item
The choice of $Q^{(L)}$ or $Q^{(R)}$, which was based on the number of leading zeros, is not guaranteed to produce better results than other choices,
although the numerical examples below seem to indicate that it performs well. 
\item
It is, in general, difficult to predict which norm provides the best result, but in many applications the size of the matrix coefficients limits
that choice to the $1$-norm or the $\infty$-norm.
\end{itemize}

We illustrate the usefulness of Theorem~\ref{Theorem_Matrix_Bakchik} and our choice of $Q^{(L)}$, defined by~(\ref{Qdef}), and compare it to Theorem~2.2  
from~\cite{Melman_RS} (the generalization to matrix polynomials of Theorem~8.3.1 in~\cite{RS}) at the hand of the following 
two examples. In the first, we generate random matrix polynomials, whereas the second one is taken from the engineering literature. 
\vskip 0.25cm

\noindent {\bf Example 1.} Here we generated $1000$ matrix polynomials with complex elements,
whose real and complex parts are uniformly randomly distributed on the interval $[-10,10]$. We then premultiplied each matrix polynomial by the inverse of its 
leading coefficient to make its leading coefficient the identity matrix.
We examined four cases with $n=20$ and $25 \times 25$ coefficients: 
$k=3,\ell=5$, $k=5,\ell=3$, $k=\ell=5$, and $k=\ell=1$, and one case with $n=4$, $250 \times 250$ coefficients, and $k=\ell=1$.
Table~\ref{Table_Example1} lists the averages of the ratios of the Cauchy radii to the modulus of the largest eigenvalue, 
i.e., the closer this number is to $1$, the better it is.
This was done for the Cauchy radius of the given matrix polynomial with the $1$-norm and five consecutive applications of Theorem~\ref{Theorem_Matrix_Bakchik}, 
labeled as level 1-5, using $Q^{(L)}$ defined by~(\ref{Qdef}) for each application. In each column, the numbers on the left are the ratios 
obtained by Theorem~\ref{Theorem_Matrix_Bakchik}, while the ones on the 
right are the ratios from~Theorem~2.2 in~\cite{Melman_RS}.
Clearly, significant improvements can obtained from Theorem~\ref{Theorem_Matrix_Bakchik}.
Moreover, the advantage of having another multiplier in addition to the one from~\cite{Melman_RS} is that it can sometimes accelerate an otherwise 
slowly progressing recursion.  
%
%
\begin{table}[H]
\begin{center}
\small           
\begin{tabular}{c|c|c|c|c|c}
Level             & $n=20$, $m=25$       & $n=20$, $m=25$       &  $n=20$, $m=25$      & $n=20$, $m=25$      & $n=4$, $m=250$            \\         
                  & $k=3$, $\ell=5$      & $k=5$, $\ell=3$      &  $k=\ell=5$          & $k=\ell=1$          & $k=\ell=1$                \\ \hline 
                  &                      &                      &                      &                     &                           \\
Cauchy            &        1.991         &        1.482         &        1.492         &       8.442         &        33.963             \\        
1                 & 1.257 $\;|\;$ 1.404  & 1.236 $\;|\;$ 1.264  & 1.165 $\;|\;$ 1.231  & 2.003 $\;|\;$ 2.880 & 3.154 $\;|\;$ 5.725       \\
2                 & 1.135 $\;|\;$ 1.198  & 1.155 $\;|\;$ 1.235  & 1.151 $\;|\;$ 1.145  & 1.419 $\;|\;$ 1.770 & 1.763 $\;|\;$ 2.419       \\
3                 & 1.127 $\;|\;$ 1.190  & 1.145 $\;|\;$ 1.217  & 1.146 $\;|\;$ 1.358  & 1.237 $\;|\;$ 1.681 & 1.361 $\;|\;$ 2.350       \\
4                 & 1.123 $\;|\;$ 1.186  & 1.117 $\;|\;$ 1.152  & 1.093 $\;|\;$ 1.130  & 1.195 $\;|\;$ 1.366 & 1.326 $\;|\;$ 1.574       \\
5                 & 1.118 $\;|\;$ 1.184  & 1.070 $\;|\;$ 1.145  & 1.087 $\;|\;$ 1.126  & 1.194 $\;|\;$ 1.328 & 1.326 $\;|\;$ 1.543       \\
\end{tabular}
\caption{Comparison of Cauchy radii for Example~1.}
\label{Table_Example1}
\end{center}
\end{table}
\normalsize

\noindent {\bf Example 2.}
%
%
This example is taken from~\cite{FPS}, where a structural dynamics model representing a reinforced concrete machine foundation is formulated as a sparse quadratic 
$3627 \times 3627$ eigenvalue problem with $k=\ell=1$. 
Of the many bounds on the eigenvalues that were examined in~\cite{HighamTisseur} for the $1$-norm and $\infty$ norm 
for this problem (the $2$-norm is too costly here), the Cauchy radius was among the best.
Theorem~2.2 in~\cite{Melman_RS} improves those bounds significantly, but Theorem~\ref{Theorem_Matrix_Bakchik} improves them even more.
Table~\ref{Table_Example2} shows the Cauchy radius and its improvements from Theorem~2.2 in~\cite{Melman_RS} and Theorem~\ref{Theorem_Matrix_Bakchik} for the 
$1$-norm on the left and the $\infty$-norm on the right. In each column, the numbers on the left are obtained from Theorem~\ref{Theorem_Matrix_Bakchik}, while 
those on the right are from Theorem~2.2 in~\cite{Melman_RS}. Here too, we have carried out five recursions of Theorem~\ref{Theorem_Matrix_Bakchik},
each time using $Q^{(L)}$ defined in~(\ref{Qdef}).
The modulus of the largest eigenvalue is $2.120 \times 10^{4}$, and in the table all bounds were divided by $10^{4}$.
%
%
\begin{table}[H]
\begin{center}
\small           
\begin{tabular}{cc}
\begin{tabular}{c|c}
Cauchy            &        3.532             \\        
Level 1           & 2.762 $\;|\;$ 3.349      \\
Level 2           & 2.427 $\;|\;$ 2.737      \\
Level 3           & 2.413 $\;|\;$ 2.722      \\
Level 4           & 2.272 $\;|\;$ 2.425      \\
Level 5           & 2.271 $\;|\;$ 2.419      \\
\end{tabular}
&
\hskip 1cm
\begin{tabular}{c|c}
Cauchy            &        3.173             \\        
Level 1           & 2.658 $\;|\;$ 3.064      \\
Level 2           & 2.380 $\;|\;$ 2.652      \\
Level 3           & 2.363 $\;|\;$ 2.598      \\
Level 4           & 2.260 $\;|\;$ 2.380      \\
Level 5           & 2.260 $\;|\;$ 2.374      \\
\end{tabular}
\\
\end{tabular}
\caption{Comparison of Cauchy radii for Example~2 with the $1$-norm (left) and the $\infty$-norm (right).}
\label{Table_Example2}
\end{center}
\end{table}
\normalsize

%
%
%
%

\section{Improved Cauchy radius for scalar polynomials}
\label{scalarpolynomials}

Since scalar polynomials are $1 \times 1$ matrix polynomials, Theorem~\ref{Theorem_Matrix_Bakchik} can be applied to them as a special case. Moreover,
because of their scalar nature, the theorem can be slightly refined, as stated in the following theorem.
%
%
\begin{theorem}
\label{Theorem_Scalar_Bakchik}
Let $p(z) = \sum_{j=0}^{n} a_{j} z^{j}$ be a polynomial of degree $n$ with complex coefficients that is at least a trinomial,
and let~$k$~and~$\ell$ be the smallest positive integers such that $\anmk$ and $\anmkmell$ are not zero. Define
\begin{eqnarray*}
& & q_{1}(z) = \lb \an z^{k+\ell} - \anmk z^{\ell} - \anmkmell \rb p(z) \; ,    \\
& & q_{2}(z) = \lb \an z^{2k} - \anmk z^{k} + \dfrac{\anmk^{2}}{\an} \rb p(z) \; ,   \\
\hskip -1cm \text{and, when $\ell=k$,} & &   \\ 
& & q_{3}(z) = \lb \an z^{2k} - \anmk z^{k} - a_{n-2k} + \dfrac{\anmk^{2}}{\an} \rb p(z) \; .     
\end{eqnarray*}
Then the following holds.
\newline {\bf (1)} For any admissible values of $k$ and $\ell$, $\rho[q_{j}] \leq \rho[p]\,$ for $j=1,2$, and if $\ell=k$, then $\rho[q_{3}] \leq \rho[p]$. 
\newline {\bf (2)} If all the coefficients of $p$ are nonzero, then the inequalities in part (1) are strict, unless $p$ has a zero of modulus $\rho[p]$.
\end{theorem}
\prf
The first part of the theorem follows immediately from Theorem~\ref{Theorem_Matrix_Bakchik} as a special case because complex numbers are $1 \times 1$ 
complex matrices.
The second part requires some elaboration. To avoid tedious repetition, we present a detailed proof only for $q_{1}$, and sketch the proof for $q_{2}$ and $q_{3}$.
Throughout, if the index of a quantity is inadmissible, then that quantity is set equal to zero.

We now assume that all the coefficients of $p$ are nonzero, so that $k=\ell=1$, and we begin with $q_{1}$.
The expression corresponding to $S(z)$, defined by~(\ref{powereq2}) in~the proof of Theorem~\ref{Theorem_Matrix_Bakchik}, is given by
\begin{eqnarray*}
S(z) & =  & -\anmo^{2}z^{n} + \an z^{2} \sum_{k=1}^{n-3}\aj z^{j} -\anmo z \sum_{j=0}^{n-2} \aj z^{j} -\anmt \sum_{j=0}^{n-1} \aj z^{j} \\
& = & -\anmo^{2}z^{n} + \sum_{j=2}^{n-1} \an\ajmt z^{j} - \sum_{j=1}^{n-1} \anmo \ajmo z^{j} - \sum_{j=0}^{n-1} \anmt \aj z^{j} \\
& = & -\anmo^{2}z^{n} + \sum_{j=2}^{n-1} \lb \an \ajmt - \anmo \ajmo - \anmt \aj \rb z^{j} - \lb \anmo \azero + \anmt \aone \rb z  - \anmt \azero \; , \\
\end{eqnarray*}
while the expression corresponding to~$\Phi$ becomes
\bdis
\Phi(z) = |\anmo|^{2}z^{n} + \sum_{j=2}^{n-1} \left | \an \ajmt - \anmo \ajmo - \anmt \aj \right | z^{j} + \left | \anmo \azero + \anmt \aone \right | z  
+ |\anmt \azero| \; . 
\edis
For $x=\rho[p]$, the inequality corresponding to~(\ref{firstsummation}) is
\beq
\label{phiineq}    
\Phi(x) \leq |\anmo |^{2}x^{n} + \sum_{j=2}^{n-1} \lb |\an \ajmt | + |\anmo \ajmo | + |\anmt \aj| \rb x^{j} + \lb |\anmo \azero | + |\anmt \aone | \rb x
+ |\anmt \azero | \; . 
\eeq
The inequality in~(\ref{phiineq}) is strict, unless 
\begin{eqnarray}
& & \hskip -0.5cm | \an \ajmt - \anmo \ajmo - \anmt \aj | = |\an \ajmt | + |\anmo \ajmo | + |\anmt \aj |  \;\; (j=2,...,n-1)  \label{abseq1}  \\
& & \hskip -1.5cm \text{and} \nonumber \\
& & \hskip -0.5cm | \anmo \azero + \anmt \aone | = |\anmo \azero | + |\anmt \aone |  \label{abseq2} \; .
\end{eqnarray}
We now define $\varphij = \text{arg} \, \aj$ and use $\varphi \cong \psi$ to indicate
that $\varphi$ and $\psi$ only differ by an integer multiple of $2\pi$, so that 
$e^{i\varphi}= e^{i\psi}$. If~(\ref{abseq1}) and~(\ref{abseq2}) hold, then we have from~(\ref{abseq1}) for $j=2,...,n-1$, that
\begin{eqnarray}
& & \varphin + \varphijmt \cong \varphinmo + \varphijmo + \pi \;, \; \text{or} \;\; 
\varphijmt \cong \varphijmo + \varphinmo - \varphin + \pi \; , 
\label{abseq3}  \\
& & \hskip -1cm \text{and}   \nonumber     \\
& & \varphinmo + \varphijmo \cong \varphinmt + \varphij \; , \label{abseq4}  
\end{eqnarray}
while from~(\ref{abseq2}) we have
\beq
\label{abseq5}
\varphinmo + \varphizero \cong \varphinmt + \varphione \; .
\eeq
Combining~(\ref{abseq3}) with the substitution $j=j-1$ in~(\ref{abseq4}), we obtain for $j=3,...,n-1$ that
\bdis
\varphijmt \cong \varphijmo + \varphinmo - \varphin + \pi \cong \varphijmo + \varphinmt - \varphinmo \; , 
\edis
which implies that 
\beq
\label{phinmteq}  
\varphinmt \cong 2\varphinmo - \varphin + \pi \; .
\eeq
Substituting this in~(\ref{abseq5}) shows that~(\ref{abseq5}) is covered by~(\ref{abseq3}).
Since~(\ref{abseq3}) is equivalent to 
\beq
\label{phieq}
\varphijmo \cong \varphij + \varphinmo - \varphin + \pi \;\;\;\; (j=1,...,n-2) \; ,
\eeq
we have obtained from~(\ref{phinmteq}) that~(\ref{phieq}) also holds for $j=n-1$.
From here on, the proof follows that of Theorem~8.3.1.in~\cite{RS}.
As in that proof, the equations in~(\ref{phieq}), used recursively for $j=n-1,...,1$, yield    
\bdis
\varphi_{n-j} \cong \varphinmo + (j-1) \Delta \;\;\;\; (j=1,...,n) \; ,
\edis
where $\Delta=\varphinmo - \varphin + \pi$, which is equivalent to
\beq    
\label{recursion}
\varphij \cong (n-j) \Delta + \varphin - \pi \;\;\;\; (j=0,...,n-1) \; .
\eeq
Using~(\ref{recursion}), we now show that, under these conditions, $xe^{i\Delta}$, where $x=\rho[p]$, is a zero of $p$:
\begin{eqnarray*}
\sum_{j=0}^{n} \aj (xe^{i\Delta})^{j} & = &  \sum_{j=0}^{n-1} |\aj|e^{i\varphij}x^{j} e^{i j \Delta} + |\an|e^{i\varphin} x^{n} e^{i n \Delta}  \\
                      & = &  e^{i(\varphin + n\Delta)} \lb \sum_{j=0}^{n-1} |\aj|e^{i(\varphij-(n-j)\Delta - \varphin)} x^{j} + |\an| x^{n} \rb  \\
                      & = &  e^{i(\varphin + n\Delta)} \lb \sum_{j=0}^{n-1} e^{-i\pi} |\aj| x^{j} + |\an| x^{n} \rb = 0 \; ,  \\
\end{eqnarray*}
and $e^{-i\pi} = -1$, so that $x=\rho[p]$ is indeed a zero of $p$.

For $q_{2}$, we obtain for $S(z)$, defined in~(\ref{powereq4}), 
\begin{eqnarray*}
S(z) & =  & \an \anmt z^{n} + \an z^{2} \sum_{k=1}^{n-3}\aj z^{j} -\anmo z \sum_{j=0}^{n-2} \aj z^{j} + \dfrac{\anmo^{2}}{\an} \sum_{j=0}^{n-1} \aj z^{j} \\
& = & \an \anmt z^{n} + \sum_{j=2}^{n-1} \an \ajmt z^{j} - \sum_{j=1}^{n-1} \anmo \ajmo z^{j} + \sum_{j=0}^{n-1} \dfrac{\anmo^{2}\aj}{\an}  z^{j} \\
& = & \an \anmt z^{n} + \sum_{j=2}^{n-1} \lb \an\ajmt - \anmo \ajmo + \dfrac{\anmo^{2}\aj}{\an} \rb z^{j} + \lb -\anmo \azero + \dfrac{\anmo^{2}\aone}{\an} \rb z
         + \dfrac{\anmo^{2}\azero}{\an}  \; , \\
\end{eqnarray*}
and for $q_{3}$, we obtain, as in~(\ref{powereq6}),
\begin{eqnarray*}
S(z) & =  & \an \anmt z^{n} + \an z^{2} \sum_{k=1}^{n-3}\aj z^{j} -\anmo z \sum_{j=0}^{n-2} \aj z^{j} 
            - \lb \anmt - \dfrac{\anmo^{2}}{\an} \rb \sum_{j=0}^{n-1} \aj z^{j} \\
& = & \an \anmt z^{n} + \sum_{j=2}^{n-1} \an \ajmt z^{j} - \sum_{j=1}^{n-1} \anmo \ajmo z^{j} - \sum_{j=0}^{n-1} \lb \anmt - \dfrac{\anmo^{2}}{\an} \rb \aj z^{j} \\
& = & \an \anmt z^{n} + \sum_{j=2}^{n-1} \lb \an \ajmt - \anmo \ajmo -\anmt \aj + \dfrac{\anmo^{2}\aj}{\an} \rb z^{j} \\
& &   \hskip 4cm + \lb -\anmo \azero - \anmt \aone  + \dfrac{\anmo^{2}\aone}{\an}  \rb z + \lb -\anmt + \dfrac{\anmo^{2}}{\an} \rb \azero  \; . \\
\end{eqnarray*}
Analogously to the proof for~$q_{1}$, we now obtain the same equations~(\ref{recursion}) for both $q_{2}$ and $q_{3}$, from which the proof follows
for these polynomials as well. \qed

Here too, and for the same reasons as in the matrix case, we make the following choice to improve the Cauchy radius of $p$:
\begin{eqnarray}
& & q(z) = \lcb 
\begin{array}{ll}
\lb \an z^{k+\ell} - \anmk z^{\ell} - \anmkmell \rb p(z)   & \text{if $\ell < k$} \; ,  \\
                                                           &                       \\
\lb \an z^{2k} - \anmk z^{k} + \anmk^{2}\an^{-1} \rb p(z)  & \text{if $\ell > k$}  \; , \\
                                                           &                       \\
\lb \an z^{2k} - \anmk z^{k} - a_{n-2k} + \dfrac{\anmk^{2}}{\an} \rb p(z)  & \text{if $\ell = k$} \; . \\
\end{array}
\right . \nonumber \\ \nonumber  
\end{eqnarray}

\noindent {\bf Acknowledgement.} We thank Fran\c coise Tisseur for sending us the matrices in Example~2.

\end{document}